# An approach to determining the existence of a limit distribution of additive arithmetic functions

Victor Volfson

ABSTRACT An approach will be proposed to determine the existence of a limit distribution of additive arithmetic functions in this work. It is based on assertions that will be proven in this work and on the properties of Dirichlet convolution and Möbius inversion. Based on this approach, formulas will be obtained for finding the mean and variance for the limit distribution of additive arithmetic functions. Examples of this approach to proving the existence of a limit distribution of additive arithmetic functions and finding the mean value and variance for the limit distribution of additive arithmetic functions are considered.

Keywords: arithmetic function, additive arithmetic function, strongly additive arithmetic function, multiplicative arithmetic function, strongly multiplicative arithmetic function, probability space, limit distribution of a real arithmetic function, mean value of an arithmetic function, variance of an arithmetic function, characteristic function of a real arithmetic function, Wintner's theorem on the mean value of a multiplicative arithmetic function, Dirichlet convolution, Möbius inversion.

# 1. INTRODUCTION

In general, an arithmetic function is a function defined on the set of natural numbers and taking values on the set of complex numbers. The name arithmetic function is due to the fact that this function expresses some arithmetic property of the natural series.

Let $m_1, m_2$ are relatively prime natural numbers.

An arithmetic function $f(m), m = 1, ..., n$ is called additive if it holds:

$$f(m_1 m_2) = f(m_1) + f(m_2), f(1) = 0. \tag{1.1}$$

An arithmetic function $g(m), m = 1, ..., n$ is called multiplicative if it holds:

$$g(m_1 m_2) = g(m_1) g(m_2), g(1) = 1. \tag{1.2}$$

Let $p$ is an arbitrary prime number and $\alpha > 1$ is a natural number.

An additive arithmetic function $f(m), m = 1, ..., n$ is strongly additive if:

$$f(p^\alpha) = f(p). \tag{1.3}$$

A multiplicative arithmetic function $g(m), m = 1, ..., n$ is strongly multiplicative if

$$g(p^\alpha) = g(p). \tag{1.4}$$

The axiomatics of probability theory makes it possible to use probability theory to find the distribution of any arithmetic functions. Any interval of a natural series $[1, n]$ can be naturally transformed into a probability space $(\Omega_n, A_n, P_n)$ by taking as $\Omega_n = \{1, 2, ..., n\}$, $A_n$ - all subsets $\Omega_n$, $P_n(A) = \#(m \in A)/n$, where $\#(m \in A)$ is the number of natural numbers in the subset $A$. Then an arbitrary arithmetic function $f(m), m = 1, ..., n$ on $\Omega_n$ can be considered as a random variable $x_n$ on this probability space:

$$x_n(m) = f(m)(1 \leq m \leq n). \tag{1.5}$$

In particular, we can talk about the average value of an arithmetic function (mathematical expectation):

$$E[x_n] = E[f,n] = \frac{1}{n}\sum_{m=1}^{n} f(m), \tag{1.6}$$

variances:

$$D[x_n] = D[f,n] = \frac{1}{n}\sum_{m=1}^{n} f^2(m) - (\frac{1}{n}\sum_{m=1}^{n} f(m))^2, \tag{1.7}$$

and about the distribution function for a real arithmetic function $f(m), m = 1,...,n$:

$$F_n(y) = P_n\{m \leq n : f(m) \leq y\} \tag{1.8}$$

and characteristic function:

$$\varphi_{x_n}(t) = \varphi_{f,n}(t) = E[e^{itf_n}] = \frac{1}{n}\sum_{m=1}^{n} e^{itf(m)}. \tag{1.9}$$

Levi's continuity theorem relates pointwise convergence characteristic functions of random variables with the convergence of these random variables values by distribution.

Direct theorem. If the sequence of distribution functions $\{F_n(y)\}$ weakly converges to the distribution function $F(y)$ as n → ∞, then the sequence of corresponding characteristic functions $\{\varphi_{f,n}(t)\}$ converges pointwise to the characteristic function $\varphi_f(t)$ as n → ∞.

Converse theorem. Let the sequence of characteristic functions $\{\varphi_{f,n}(t)\}$ converge pointwise to a function $\varphi_f(t)$ continuous at point 0. Then the sequence of corresponding distribution functions $\{F_n(y)\}$ weakly converges to the function $F(y)$ and $\varphi_f(t)$ is a characteristic function corresponding to the distribution function $F(y)$.

Erdős and Wintner proved [1] that the condition for the convergence of series for a positive value $R$ for a real additive arithmetic function $f(m), m = 1,...,n$ is:

$$\sum_{|f(p)| \leq R} \frac{f(p)}{p},$$

$$\sum_{|f(p)| \leq R} \frac{f^2(p)}{p},$$

$$\sum_{|f(p)|>R} \frac{1}{p}.$$

is sufficient for the convergence of the sequence of distribution functions $\{F_n(y)\}$ as n → ∞ to the limit distribution function $F(y)$.

The idea of the proof is based on theorems on sums of independent random variables. It should be noted that if the indicated series converge at some positive value $R$, then without loss of generality we can take $R=1$.

When $|f(p)| \leq 1$ then it is executed:

$$\sum_{p} \frac{1}{p} < \infty,$$

$$\sum_{p} \frac{f^2(p)}{p} \leq \sum_{p} \frac{|f(p)|}{p} < \infty$$

and the conditions of the Erdős-Wintner theorem are simplified - only convergence of the series $\sum_{p} \frac{|f(p)|}{p}$ is required.

It is also interesting to prove the fact that the distribution function converges to the limit distribution function and to find the asymptotic distribution laws for additive arithmetic functions in the case when the indicated series diverge. Using the central limit theorem, Erdős and Kac [2] found that the limit distribution law for a strongly additive arithmetic function $\omega(m), m=1,...,n$ is normal. Later they proved that this holds for all strongly additive arithmetic functions for which $B(n) = \sqrt{\sum_{p \leq n} \frac{f^2(p)}{p}} \to \infty$ as n → ∞ and $|f(p)| \leq 1$.

Kubilius [3] introduced the class H of additive arithmetic functions for which the limit distribution of additive arithmetic functions coincides with the limit distribution of strongly additive arithmetic functions. Mean values, variances and moments of higher orders were found for this class of additive arithmetic functions in [4.5].

An approach to determining the existence of a limit distribution of additive arithmetic functions will be proposed in this work, different from those discussed above. This approach based on assertions that will be proven in the second chapter of this work and on the properties of the Dirichlet convolution and Möbius inversion indicated in the third chapter of the work. The application of the properties of Dirichlet convolutions and Möbius inversion to the proof of the existence of the limit distribution of additive arithmetic functions will be given in the fourth chapter of the work. Examples of this approach to proving the existence of the limit distribution of additive arithmetic functions will be considered in the fifth chapter. Based on this approach, formulas for finding the mean value and variance for the limit distribution of additive arithmetic functions will be obtained in the sixth chapter of the work.

2. AN APPROACH TO DETERMINING THE EXISTENCE OF A LIMIT DISTRIBUTION ADDITIVE ARITHMETIC FUNCTIONS

Theorems on the mean value of a multiplicative function provide an approach not only to analytic number theory, but also to probabilistic number theory.

This approach was applied and proof of the existence of the limit distribution of a real additive arithmetic function was given on the mean value of a multiplicative arithmetic function using Delange's theorem in the work [6]. However, this proof is quite complicated. A simpler and, it seems to me, more elegant implementation of this approach will be shown here – through the well-known Wintner theorem [7].

Let us prove the following assertion using Assertion 1.

Assertion 1

Let $g$ is a complex-valued multiplicative arithmetic function for which the following condition is satisfied:

$$\sum_{n=1}^{\infty} \frac{|g(n) * \mu(n)|}{n} < \infty, \qquad (2.1)$$

where * is the Dirichlet convolution. Then the limit of the average value is:

$$\lim_{x \to \infty} \frac{1}{x} \sum_{n \leq x} g(n) = \prod_p (1 - \frac{1}{p}) \sum_{\alpha \geq 0} \frac{g(p^{\alpha})}{p^{\alpha}}. \qquad (2.2)$$

Assertion 2

Let $f(n)$ is a real additive arithmetic function for which the following condition is

satisfied:

$$\sum_{n=1}^{\infty} \frac{|e^{itf(n)} * \mu(n)|}{n} < \infty. \quad (2.3)$$

Then $f(n)$ has a limit distribution as n $\to \infty$, for which the characteristic function has the form:

$$\varphi_f(t) = \prod_p (1 - \frac{1}{p}) \sum_{\alpha \geq 0} \frac{e^{itf(p^\alpha)}}{p^\alpha}. \quad (2.4)$$

Proof

Using the probability space $(\Omega_n, A_n, P_n)$ specified in Chapter 1, based on (1.5), a real additive arithmetic function $f(n)$ can be considered as a random variable for which, in accordance with (1.8), there is a distribution function $F_n$ and, in accordance with (1.9), a characteristic function - $\varphi_{f,n}(t)$.

If $f(n)$ is a real additive arithmetic function, then $e^{itf(n)}$ is a complex-valued multiplicative arithmetic function and if $e^{itf(n)}$ satisfies condition (2.3), then condition (2.1) of Assertion 1 is satisfied. Then, based on Assertion 1, (2.2) is satisfied.

Let us substitute a complex-valued multiplicative arithmetic function $e^{itf(n)}$ into (2.2) and obtain:

$$\lim_{x \to \infty} \frac{1}{x} \sum_{n \leq x} e^{itf(n)} = \prod_p (1 - \frac{1}{p}) \sum_{\alpha \geq 0} \frac{e^{itf(p^\alpha)}}{p^\alpha}. \quad (2.5)$$

We find the limit of the characteristic function of a real additive arithmetic function $f(n)$ on the left in (2.5), and on the right - the function $\varphi_f(t)$, which is continuous and at $t = 0$ is equal to:

$$\prod_p (1 - \frac{1}{p}) \sum_{\alpha \geq 0} \frac{1}{p^\alpha} = \prod_p (\frac{p-1}{p})(\frac{p}{p-1}) = 1. \quad (2.6)$$

Thus, based on the inverse Lévy theorem, $\varphi_f(t)$ is the characteristic function for the limit distribution $F$ of a real additive arithmetic function $f(n)$ as n $\to \infty$, which corresponds to Assertion 2.

Let us consider condition (2.3) of Assertion 2 in more detail.

It is known that the condition for the absolute convergence of Dirichlet series of any multiplicative arithmetic functions $g_1, g_2$ is sufficient for the absolute convergence of the Dirichlet series of $g_1 * g_2$. However, the Dirichlet series $\sum_{n=1}^{\infty} \frac{\mu(n)}{n}$ converges only conditionally and the Dirichlet series $\sum_{n=1}^{\infty} \frac{e^{itf(n)}}{n}$ also does not converge absolutely. Therefore, to check the fulfillment of property (2.3), it is necessary to analyze the Möbius inversion $e^{itf(n)} * \mu(n)$ for a specific real additive arithmetic function $f$.

Arithmetic functions $e^{itf(n)}$ and $\mu(n)$ are multiplicative. Therefore, it is sufficient to consider the convolution of multiplicative arithmetic functions on the product of powers of prime numbers. This analysis requires some properties of Dirichlet convolution and Möbius inversion, which we will look at in the next chapter.

## 3. PROPERTIES OF DIRICHLET CONVOLUTION AND MOBIUS INVERSION FOR MULTIPLICATIVE ARITHMETIC FUNCTIONS

Let us consider the indicated properties for arbitrary multiplicative arithmetic functions.

Assertion 3

The Dirichlet convolution for multiplicative arithmetic functions $h, g$ and an arbitrary prime $p$ is equal to:

$$h * g(p) = g(p) + h(p). \tag{3.1}$$

Möbius inversion for a multiplicative arithmetic function $h$ and an arbitrary prime $p$ is equal to:

$$h * \mu(p) = h(p) - 1. \tag{3.2}$$

The Dirichlet convolution for multiplicative arithmetic functions $h, g$ and $n \geq 2$ is equal to:

$$h * g(p^n) = g(p^n) + h(p)g(p^{n-1}) + ... + h(p^n). \tag{3.3}$$

Möbius inversion for a multiplicative arithmetic function $h$ and $n \geq 2$ equals:

$$h * \mu(p^n) = h(p^n) - h(p^{n-1}). \tag{3.4}$$

When $n = p_1^{a_1}...p_k^{a_k}$ Dirichlet convolution for multiplicative arithmetic functions $h, g$ is equal to:

$$h * g(p_1^{a_1}...p_k^{a_k}) = \prod_{l=1}^{k}(g(p_l^{a_l}) + h(p_l)g(p_l^{a_l-1}) + ... + h(p_l^{a_l})). \tag{3.5}$$

For $n = p_1^{a_1}...p_k^{a_k}$ the Möbius inversion for a multiplicative arithmetic function $h$ is:

$$h * \mu(p_1^{a_1}...p_k^{a_k}) = \prod_{l=1}^{k}(h(p_l^{a_l}) - h(p_l^{a_l-1})). \tag{3.6}$$

Proof

It is known the formula for Dirichlet convolution:

$$h * g(n) = \sum_{d/n} h(d)g(\frac{n}{d}). \tag{3.7}$$

Based on (3.7) we have for multiplicative arithmetic functions $h, g$ and an arbitrary prime $p$:

$$h * g(p) = h(1)g(p) + h(p)g(1) = g(p) + h(p). \tag{3.8}$$

Taking into account (3.8) we obtain the Möbius inversion for a multiplicative arithmetic function $h$ and an arbitrary prime $p$:

$$h * \mu(p) = h(p) - 1. \tag{3.9}$$

Based on (3.7) we have for multiplicative arithmetic functions $h, g$ and $n \geq 2$:

$$h * g(p^n) = h(1)g(p^n) + h(p)g(p^{n-1}) + ... + h(p^n)g(1) = g(p^n) + h(p)g(p^{n-1}) + ... + h(p^n). \tag{3.10}$$

Having in mind (3.10) we obtain the Möbius inversion for the multiplicative arithmetic function $h$ and $n \geq 2$:

$$h * \mu(p^n) = h(p^n) - h(p^{n-1}). \tag{3.11}$$

When $n = p_1^{a_1}...p_k^{a_k}$, based on (3.10), we have for multiplicative arithmetic

functions $h, g$ :

$$h * g(p_1^{a_1}...p_k^{a_k}) = \prod_{l=1}^{k}(g(p_l^{a_l}) + h(p_l)g(p_l^{a_l-1}) + ... + h(p_l^{a_l})),$$

which corresponds to (3.5).

For $n = p_1^{a_1}...p_k^{a_k}$ based on (3.11) we obtain the Möbius inversion for multiplicative arithmetic function $h$:

$$h * \mu(p_1^{a_1}...p_k^{a_k}) = \prod_{l=1}^{k}(h(p_l^{a_l}) - h(p_l^{a_l-1})),$$

which corresponds to (3.6).

If the arithmetic function $h$ is strongly multiplicative and $n \geq 2$, then based on (3.3) and (3.4) we have:

$$h * g(p^n) = g(p^n) + h(p)(g(p^{n-1}) + g(p^{n-2}) + ... + g(p) + 1), \tag{3.12}$$

and

$$h * \mu(p^n) = h(p^n) - h(p^{n-1}) = h(p) - h(p) = 0. \tag{3.13}$$

4. APPLICATION OF THE PROPERTIES OF DIRICHLET CONVOLUTION AND MOBIUS INVERSION TO THE PROOF OF THE EXISTENCE OF THE LIMIT DISTRIBUTION OF ADDITIVE ARITHMETIC FUNCTIONS

Based on (3.4), Assertion 2 can be written in the form.

Assertion 4

Let $f(n)$ is a real additive arithmetic function for which the following condition is satisfied:

$$\sum_{p^\alpha} \frac{|e^{itf(p^\alpha)} - e^{itf(p^{\alpha-1})}|}{p^\alpha} < \infty. \tag{4.1}$$

Then $f(n)$ has a limit distribution as n $\to \infty$, for which the characteristic function has the form:

$$\varphi_f(t) = \prod_p (1-\frac{1}{p}) \sum_{\alpha \geq 0} \frac{e^{itf(p^\alpha)}}{p^\alpha}. \qquad (4.2)$$

Assertion 5

Let $f^*(n)$ is a real strongly additive arithmetic function for which the following conditions are satisfied for prime values $p$:

$$\sum_p \frac{|e^{itf^*(p)} - 1|}{p} < \infty. \qquad (4.3)$$

Then $f^*(n)$ has a limit distribution as $n \to \infty$, for which the characteristic function has the form:

$$\varphi_{f^*}(t) = \prod_p (1-\frac{1}{p})(1 + \frac{e^{itf^*(p)}}{p}). \qquad (4.4)$$

Proof

Based on the Möbius inversion property (3.2), the following holds for an arbitrary prime $p$:

$$e^{itf^*(p)} * \mu(p) = e^{itf^*(p)} - 1. \qquad (4.5)$$

Taking into account (3.13) it holds for an arbitrary prime $p$ and $\alpha \geq 2$:

$$e^{itf^*(p^\alpha)} * \mu(p^\alpha) = 0. \qquad (4.6)$$

Based on (4.5) and (4.6), Assertion 4 is satisfied, and therefore (4.2) is true.

Taking into account (4.2), and that $f^*(n)$ is a real strongly additive arithmetic function, we obtain (4.4).

Assertion 6

Let $f^*(n)$ is a real strongly additive arithmetic function for which $|f^*(p)| \leq 1$ and the following condition is satisfied:

$$\sum_p \frac{|f^*(p)|}{p} < \infty. \qquad (4.7)$$

Then $f^*(n)$ has a limit distribution as $n \to \infty$, for which the characteristic function has the form:

$$\varphi_{f^*}(t) = \prod_p (1-\frac{1}{p})(1+\frac{e^{itf^*(p)}}{p}). \tag{4.8}$$

Proof

We use that $|e^{-itf^*(p)/2}|=1$, $|f^*(p)|\leq 1$, and multiply:

$$|e^{itf^*(p)} - 1 \| e^{-itf^*(p)/2}| = |e^{itf^*(p)/2} - e^{-itf^*(p)/2}| = |2i\sin(tf^*(p)/2)| \leq |tf^*(p)|. \tag{4.9}$$

Based on (4.9) we obtain:

$$\sum_p \frac{|e^{itf^*(p)} - 1|}{p} \leq \sum_p \frac{|tf^*(p)|}{p} < \infty. \tag{4.10}$$

Inequality (4.10) is satisfied if

$$\sum_p \frac{|f(p)|}{p} < \infty. \tag{4.11}$$

Therefore, in this case, condition (4.3) of Assertion 5 is satisfied, which entails the fulfillment of (4.8).

Assertion 7

Let $f(n)$ is a real additive arithmetic function for which $|f^*(p)|\leq 1$ and the following condition is satisfied:

$$\sum_p \frac{|f(p)|}{p} < \infty. \tag{4.12}$$

Then $f(n)$ has a limit distribution as n → ∞, for which the characteristic function has the form:

$$\varphi_f(t) = \prod_p (1-\frac{1}{p}) \sum_{\alpha \geq 0} \frac{e^{itf(p^\alpha)}}{p^\alpha}. \tag{4.13}$$

Proof

We use that $|e^{-itf(p^{k+1/2})}|=1$, $|f^*(p)|\leq 1$, and divide:

$$\frac{|e^{itf(p^k)} - e^{itf(p^{k-1})}|}{|e^{-itf(p^{k+1/2})}|} = |e^{itf(p/2)} - e^{-itf(p/2)}| = 2|\sin(itf(p)/2)| \leq |tf(p)|. \tag{4.14}$$

Based on (4.14) we obtain:

$$\sum_p \frac{|e^{itf(p^k)} - e^{itf(p^{k-1})}|}{p^k} \leq \sum_p \frac{|tf(p)|}{p^k} < \infty. \qquad (4.15)$$

Inequality (4.15) is satisfied if

$$\sum_p \frac{|f(p)|}{p^k} < \infty. \qquad (4.16)$$

If condition (4.16) is satisfied for $k=1$, then it is satisfied for $k>1$, so it is sufficient to satisfy condition (4.12). Consequently, in this case, condition (4.1) of Assertion 4 is satisfied, which entails the fulfillment of (4.13).

5. EXAMPLES OF THIS APROACH TO PROVING THE EXISTENCE OF A LIMIT DISTRIBUTION OF ADDITIVE ARITHMETIC FUNCTIONS

Let's consider an example: $f^*(n) = \ln \frac{\varphi(n)}{n}$ - a real strongly additive arithmetic function.

Value $|f^*(p)| = |\ln \frac{\varphi(p)}{p}| = |\ln(1 - \frac{1}{p})| < 1$, and series:

$$\sum_p \frac{|f^*(p)|}{p} = \sum_p \frac{|\ln(1 - \frac{1}{p})|}{p} \sim \sum_p \frac{1}{p^2} < \infty.$$

Thus, the conditions of Assertion 6 are satisfied, and therefore $f^*(n) = \ln \frac{\varphi(n)}{n}$ has a limit distribution as $n \to \infty$, for which the characteristic function has the form:

$$\varphi_{f^*}(t) = \prod_p (1 - \frac{1}{p})(1 + \frac{e^{itf^*(p)}}{p}) = \prod_p (1 - \frac{1}{p})(1 + \frac{e^{it\ln(1-1/p)}}{p}).$$

Let's consider another example: $f(n) = \ln \frac{\sigma(n)}{n}$ - a real additive arithmetic function.

Value $|f(p)| = |\ln \frac{\sigma(p)}{p}| = |\ln(1 + \frac{1}{p})| < 1$, and series:

$$\sum_p \frac{|f(p)|}{p} = \sum_p \frac{|\ln(1 + \frac{1}{p})|}{p} \sim \sum_p \frac{1}{p^2} < \infty.$$

Thus, the conditions of Assertion 7 are satisfied, and therefore $f(n) = \ln \dfrac{\sigma(n)}{n}$ has a limit distribution as $n \to \infty$, for which the characteristic function has the form:

$$\varphi_f(t) = \prod_p (1-\frac{1}{p}) \sum_{\alpha \geq 0} \frac{e^{itf(p^\alpha)}}{p^\alpha} = \prod_p (1-\frac{1}{p}) \sum_{\alpha \geq 0} \frac{e^{it\sigma(p^\alpha)}}{p^\alpha} = \prod_p (1-\frac{1}{p}) \sum_{\alpha \geq 0} \frac{e^{it(1+p+...+p^\alpha)}}{p^\alpha}.$$

## 6. MEANS AND VARIANCES OF THE LIMIT DISTRIBUTION OF ADDITIVE ARITHMETIC FUNCTIONS

Assertion 8

Let $f(n)$ is a real additive arithmetic function for which the following condition is satisfied:

$$\sum_{n=1}^{\infty} \frac{|e^{itf(n)} * \mu(n)|}{n} < \infty,$$

then the average value $f(n)$ as $n \to \infty$ is:

$$E[f] = \sum_{p^\alpha, \alpha \geq 1} \frac{f(p^\alpha)}{p^\alpha}(1-\frac{1}{p}), \qquad (6.1)$$

and the variance $f(n)$ as $n \to \infty$ is equal to:

$$D[f] = \sum_{p^\alpha, \alpha \geq 1} \frac{f^2(p^\alpha)}{p^\alpha}(1-\frac{1}{p}) - (\sum_{p^\alpha, \alpha \geq 1} \frac{f(p^\alpha)}{p^\alpha}(1-\frac{1}{p}))^2. \qquad (6.2)$$

Proof

Characteristic function of the limit distribution $f(n)$ as $n \to \infty$ (if Assertion 2 is satisfied) has the form:

$$\varphi_f(t) = \prod_p (1-\frac{1}{p}) \sum_{\alpha \geq 0} \frac{e^{itf(p^\alpha)}}{p^\alpha}, \qquad (6.3)$$

and it is the product of characteristic functions of random variables:

$$\varphi_p(t) = (1-\frac{1}{p}) \sum_{\alpha \geq 0} \frac{e^{itf(p^\alpha)}}{p^\alpha} = \sum_{\alpha \geq 0} \frac{(1-\frac{1}{p})e^{itf(p^\alpha)}}{p^\alpha}. \qquad (6.4)$$

Each such random variable is discrete, so its characteristic function can be represented as:

$$\varphi_p(t) = \sum_{\alpha=1}^{\infty} e^{itf(p^\alpha)} P_a, \qquad (6.5)$$

where $P_a$ is the probability that the random variable takes the value $e^{itf(p^\alpha)}$.

Comparing (6.4) and (6.5) we find that the random variable $f_{p_\alpha}$ takes a value $f(p^\alpha)$ with probability - $P_a = \dfrac{(1-1/p)}{p^\alpha}$. Therefore, the average value of this random variable is:

$$E[f_{p_\alpha}] = \sum_{\alpha \geq 1} \frac{f(p^\alpha)}{p^\alpha}(1 - \frac{1}{p}), \qquad (6.6)$$

and the variance of the random variable $f_{p_\alpha}$ in this case is equal to:

$$D[f_{p_\alpha}] = \sum_{\alpha \geq 1} \frac{f^2(p^\alpha)}{p^\alpha}(1 - \frac{1}{p}) - (\sum_{\alpha \geq 1} \frac{f(p^\alpha)}{p^\alpha}(1 - \frac{1}{p}))^2. \qquad (6.7)$$

Based on (6.3), when Assertion 2 is satisfied, the additive arithmetic function $f(n)$ as $n \to \infty$ is the sum of independent random variables $f_p$, therefore, using (6.6), we obtain that the average value $f(n)$ as $n \to \infty$ is equal to:

$$E[f] = \sum_{p^\alpha, \alpha \geq 1} \frac{f(p^\alpha)}{p^\alpha}(1 - \frac{1}{p}), \text{ which corresponds to (6.1),}$$

and using (6.5) we find that the variance $f(n)$ as $n \to \infty$ is equal to:

$$D[f] = \sum_{p^\alpha, \alpha \geq 1} \frac{f^2(p^\alpha)}{p^\alpha}(1 - \frac{1}{p}) - (\sum_{p^\alpha, \alpha \geq 1} \frac{f(p^\alpha)}{p^\alpha}(1 - \frac{1}{p}))^2, \text{ which corresponds to (6.2).}$$

Considering that the condition is satisfied for a strongly additive arithmetic function $f^*(n)$ - $f^*(p^\alpha) = f^*(p)$, then, based on (6.1), the average value of a strongly additive arithmetic function $f^*(n)$ as $n \to \infty$ is equal to:

$$E[f^*] = \sum_p \frac{f^*(p)}{p}(1 - \frac{1}{p}), \qquad (6.8)$$

and based on (6.2), the variance $f^*(n)$ as $n \to \infty$ is equal to:

$$D[f^*] = \sum_p \frac{f^{*2}(p)}{p}(1 - \frac{1}{p}) - (\sum_p \frac{f^*(p)}{p}(1 - \frac{1}{p}))^2. \qquad (6.9)$$

Let's look at an example of using assertion 8.

Let it be necessary to determine the mean and variance of a strongly additive arithmetic function $f^*(n) = \ln \frac{\varphi(n)}{n}$ as n → ∞.

The condition of assertion 8 is satisfied for a real, strongly additive arithmetic function $f^*(n) = \ln \frac{\varphi(n)}{n}$ (see Chapter 4). Therefore, based on (6.8), we find the average value $f^*(n) = \ln \frac{\varphi(n)}{n}$ as n → ∞:

$$E[f^*] = \sum_p \frac{f^*(p)}{p}(1-\frac{1}{p}) = \sum_p \frac{\ln(\frac{\varphi(p)}{p})}{p}(1-\frac{1}{p}) = \sum_p \frac{\ln(1-\frac{1}{p})}{p}(1-\frac{1}{p}) \sim -\sum_p \frac{1}{p^2} \text{ - converges.}$$

Based on (6.9), we find the variance $f^*(n) = \ln \frac{\varphi(n)}{n}$ as n → ∞:

$$D[f^*] = \sum_p \frac{f^{*2}(p)}{p}(1-\frac{1}{p}) - (\sum_p \frac{f^*(p)}{p}(1-\frac{1}{p}))^2 =$$

$$= \sum_p \frac{\ln^2(1-\frac{1}{p})}{p}(1-\frac{1}{p}) - (\sum_p \frac{\ln(1-\frac{1}{p})}{p}(1-\frac{1}{p}))^2 \sim \sum_p \frac{1}{p^3} - \sum_p \frac{1}{p^4} \text{ - converges.}$$

## 7. CONCLUSION AND SUGGESTIONS FOR FURTHER WORK

The next article will continue to study the asymptotic behavior of some arithmetic functions.

## 8. ACKNOWLEDGEMENTS

Thanks to everyone who has contributed to the discussion of this paper. The author is very pleased that Professor Gérald Tenenbaum showed interest in this article and read it in full before publication in the Archive. I am grateful to everyone who expressed their suggestions and comments in the course of this work.

References

1. P.Erdos, A. Wintner, Additive arithmetical functions and statistical independence. Amer.J.Math.,1939, 61. 713-721.

2. P.Erdos, M Kac. The Gausian law of errors in the theorybnof additive number-theoretic functions. Proc. Nat. Acad. Sci. U.S.A. 1939, 25, 206-207.

3. Kubilius I.P. Probabilistic Methods in Number Theory, Vilnius, Gospolitnauchizdat, Lithuanian SSR, 1962, 220 pp.

4. Volfson V.L. Estimate of the asymptotic behavior of the moments of arithmetic functions having limiting normal distribution , arXiv preprint https://arxiv.org/abs/2104.10164 (2021)

5. Volfson V.L. An estimate of asymptotics of the moments of additive arithmetic functions with a limit distribution defined on a subset of the natural series, arXiv preprint https://arxiv.org/abs/2106.10237 (2021)

6. Gérald Tenenbaum. Introduction to analytic and probabilistic number theory. Third Edition, Translated by Patrick D. F. Ion, 2015 - 630 p.

7. E. Wirsing, "Das asymptotische Verhalten von Summen über multiplikative Funktionen, II" Acta Math. Acad. Sci. Hung. , 18 (1967) pp. 411–467